\documentclass{amsart}
\usepackage{amssymb}
\newtheorem{theorem}{Theorem}[section]
\newtheorem{lemma}[theorem]{Lemma}
\newtheorem{remark}[theorem]{Remark}
\newtheorem{proposition}[theorem]{Proposition}
\newtheorem{corollary}[theorem]{Corollary}
\theoremstyle{definition}

\newtheorem{example}[theorem]{Example}
\numberwithin{equation}{section}

\def\IC{{\mathbb C}}
\def\IR{{\mathbb R}}
\def\cL{{\mathcal L}}
\def\cV{{\mathcal V}}
\def\cF{{\mathcal F}}
\def\cH{{\mathcal H}}
\def\cE{{\mathcal E}}
\def\cH{{\mathcal H}}
\def\cP{{\mathcal P}}
\def\cK{{\mathcal K}}
\def\cR{{\mathcal R}}
\def\cS{{\mathcal S}}
\def\bV{{\bf V}}
\def\diag{{\rm diag}\,}
\def\conv{{\rm conv}\,}
\def\rank{{\rm rank}\,}
\def\cB{{\mathcal B}}
\def\BH{{\cB(\cH)}}
\def\Re{{\rm Re}\,}
\def\re{{\rm Re}\,}
\def\Im{{\rm Im}\,}
\def\im{{\rm Im}\,}
\def\Ra{{\Rightarrow}}
\def\[{\left [}
\def\]{\right ]}
\def\({\left (}
\def\){\right )}
\def\la{{\langle}}
\def\ra{{\rangle}}
\def\Ra{{\ \Rightarrow\ }}

\def\cl{{\bf Cl}}
\def\int{{\bf Int}}
\def\ker{{\rm ker\ }}

\begin{document}
\openup 1\jot

\title[Higher rank numerical ranges]
{Higher rank numerical ranges and low \\
rank perturbations of quantum channels}

\author{Chi-Kwong Li}
\address{Department of Mathematics, College of William \& Mary,
Williamsburg, VA 23185}
\email{ckli@math.wm.edu}

\author{Yiu-Tung Poon}
\address{Department of Mathematics, Iowa State University,
Ames, IA 50051}
\email{ytpoon@iastate.edu}

\author{Nung-Sing Sze}
\address{Department of Mathematics, University of Connecticut,
Storrs, CT 06269}
\email{sze@math.uconn.edu}
\thanks{}


%\date{}

\dedicatory{}

%    "Communicated by" -- provide editor's name; required.
\commby{}

\begin{abstract}
For a positive integer $k$,
the rank-$k$ numerical range $\Lambda_k(A)$ of an operator $A$ acting
on a Hilbert space $\cH$ of dimension at least $k$
is the set of scalars $\lambda$ such
that $PAP = \lambda P$ for some rank $k$ orthogonal projection $P$.
In this paper, a close connection
between  low rank perturbation of an operator $A$ and $\Lambda_k(A)$
is established.
In particular, for $1 \le r < k$
it is shown that $\Lambda_k(A) \subseteq \Lambda_{k-r}(A+F)$
for any operator $F$ with $\rank (F) \le r$.
In quantum computing, this result implies that
a quantum channel  with a  $k$-dimensional error correcting
code under a perturbation of rank  $\le r$
 will still have a $(k-r)$-dimensional error correcting
code. Moreover, it is shown that if $A$ is normal or if the dimension
of $A$ is finite, then $\Lambda_k(A)$  can be obtained as the intersection
of $\Lambda_{k-r}(A+F)$ for a collection of rank $r$ operators $F$.
Examples are given to show that the result fails if $A$ is
a general operator.  The closure and the interior
of the convex set $\Lambda_k(A)$ are
completely determined. Analogous results are obtained  for
$\Lambda_\infty(A)$ defined as the set of scalars $\lambda$ such that
$PAP = \lambda P$ for an infinite rank orthogonal projection $P$.
It is shown that $\Lambda_\infty(A)$ is the intersection
of all $\Lambda_k(A)$ for $k = 1, 2, \dots$.
If $A - \mu I$ is not compact for any $\mu \in \IC$, then
the closure and the interior of $\Lambda_\infty(A)$ coincide with
those of the essential numerical range of $A$.
The situation for the special case when $A-\mu I$ is compact
for some $\mu \in \IC$ is also studied.
\end{abstract}

\maketitle

\section{Introduction}
\setcounter{equation}{0}

Let $\BH$ be the algebra of bounded linear operators acting on
a Hilbert space $\cH$. We identify $\BH$ with $M_n$ if
$\cH$ has dimension $n$.
For %a positive integer
$k \le {\rm dim}\,\cH$ , define the
{\it rank-$k$ numerical range} of $A \in \BH$ by
$$\Lambda_k(A) = \{ \lambda \in \IC: PAP = \lambda P \hbox{\rm
~~ for some rank-$k$ orthogonal projection } P \in \BH \}.$$
Note that we allow $k=\infty$ if dim $\cH=\infty$.
Evidently, $\lambda \in \Lambda_k(A)$ if and only if there is an
orthogonal basis of $\cH$ such that $\lambda I_k$
is the leading principal  submatrix of the operator matrix of $A$
with respect to the basis; equivalently, there is an
isometry $X: \IC^k \rightarrow \cH$ such that $X^*AX = \lambda I_k$.
(For $k=\infty$, we take $X:\ell^2\rightarrow \cH$.) When $k = 1$, this
concept reduces to the {\it classical numerical range} of $A$ defined by
$$W(A) = \{\la Ax,x\ra: x \in \cH, \ \la x,x\ra = 1\},$$
which is useful in studying operators and matrices; for example
see \cite{GR}.

The higher rank numerical range was introduced in connection to
the construction of quantum error correction code
in the study of quantum information theory; see \cite{Cet2}.
In quantum computing, information is stored in qubits (quantum bits).
Mathematically, the state of a qubit is represented by 
a $2\times 2$ rank one Hermitian matrix $Q$ satisfying $Q^2 = Q$.
A state of $N$-qubits $Q_1, \dots, Q_N$ is represented by their
tensor products in $M_n$ with $n = 2^N$.  A quantum channel for states of
$N$-qubits corresponds to trace preserving completely positive linear
map $\Phi: M_n \rightarrow M_n$. By the structure theory of
completely positive linear map \cite{C},  there are $T_1, \dots, T_m \in M_n$
with $\sum_{j=1}^m   T_j^*T_j = I_n$ such that
\begin{equation}\label{channel}
\Phi(X) = \sum_{j=1}^m T_jXT_j^*.
\end{equation}
Let $\bV$ be a subspace of $\IC^n$ and $P_{\bV}$ the orthogonal projection of $\IC^n$ onto $\bV$.
Then $\bV$ is a quantum error correction code for $\Phi$ if there exists a trace preserving completely positive linear map $\Psi: M_n \rightarrow M_n$ such that  $\Psi\circ\Phi(A)=A$ for all $A\in P_{\bV}M_n P_{\bV}$.
This happens if and only if there are scalars $\gamma_{ij}$
with $1 \le i, j \le m$  such that
$$P_{\bV}T_i^*T_jP_{\bV} = \gamma_{ij}P_{\bV}, \quad 1 \le i, j \le m;$$
 see \cite{Cet2,KL}. It turns out that even for a single matrix
$A$, determining $\Lambda_k(A)$ is highly non-trivial,
and the results are useful in quantum computing, say,
in constructing binary unitary channels; see \cite{Cet0}.
In a sequence of papers \cite{Cet,Cet0,Cet1,Cet2,GLW,LPS,LS,W1},
researchers studied the set $\Lambda_k(A)$ for $A \in \BH$.
Many interesting results (see P1-P8 below) were obtained. 
%for finite dimensional operators.

In the study of operator theory and applications, it is often useful to study
the properties of an operator which are stable under different kinds of
perturbation. For example,  the {\it essential numerical range} of
an infinite dimensional operator $A \in \BH$ can be defined as
\begin{equation}\label{ess}
W_e(A)=\cap\{\cl(W(A+F)): F \in \BH \hbox{ has
finite rank} \},\end{equation}
which captures many important properties
of $A$ (see  \cite{AS,FSW,SW,W2}).
Here $\cl(S)$ denotes the closure of the set $S$.
In fact, one can include all compact operators $F$ in $\BH$
on the right hand side of (\ref{ess}). If $\cK$ is the
algebra of compact operators in $\BH$
and if  $\psi:\BH\mapsto \BH/\cK$ is the canonical homomorphism of
$\BH$ onto the Calkin algebra $\BH/\cK$,
then $W_e(A)$ is the closure of
the numerical range of $\psi(A)$.
In \cite[Theorem 4]{AS}, it was also proven that
\begin{equation}\label{Linf}
\Lambda_\infty(A)=\cap\{W(A+F): F \in \BH \hbox{
has finite rank} \}.\end{equation}

In this paper, we study the change of
the higher rank numerical range of an operator
under low rank perturbation. For instance,
we show in Theorem \ref{inclusion} that
for $1 \le r < k <\infty$, if
$A, F \in \BH$ with $\rank (F) \le r$, then
\begin{equation}\label{key}
\Lambda_k(A) \subseteq \Lambda_{k-r}(A+F).
\end{equation}
In Theorem \ref{inf},
we refine the set equalities  (\ref{ess}) and (\ref{Linf})
by using  smaller sets of operators $F$ for the intersection
on the right hand sides of the equalities.

It is worth noting that the inclusion (\ref{key}) has the following implication
in the theory of quantum computing. Suppose $A \in M_n$ corresponds to a
quantum channel with a $k$-dimensional error correcting code (realized as a
subspace of $\IC^n$), then for any perturbation of the channel
$A$ by an operator $F$ of rank bounded by $r$,
the resulting channel $A+F$ will have a $(k-r)$-dimensional
error correcting code. More generally, if the matrices
$T_1, \dots, T_m$ correspond to quantum channel (\ref{channel})
with a $k$-dimensional error correcting code, and if
$T_{j}$ is changed to $T_j + F_j$ such that
the sum of the range spaces of
$$(T_i+F_i)^*(T_j+F_j)-T_i^*T_j = T_i^*F_j + F_iT_j^* + F_iF_j^*, \qquad
1 \le i, j \le n,$$
has dimension bounded by $r$, then the resulting quantum channel
will still have a $(k-r)$-dimensional error correcting code.

Our paper is organized as follows.  First, we study $\Lambda_k(A)$
for $A \in \BH$ when $k$ is
finite in Sections 2 -- 4. In Section 2, we give a complete description of
the closure and interior of $\Lambda_k(A)$.
In Section 3, we  establish inclusion (\ref{key})
for any operators $A, F \in \BH$ with rank$(F) \le r$,
where $1\le r<k < \infty$. It follows that
\begin{equation}\label{per}\Lambda_k(A) \subseteq
\cap\{\Lambda_{k-r}(A+F):F\in \BH\mbox{ has rank }\le r\}\,.
\end{equation}
In particular,  taking $r = k-1$, we have
\begin{equation}
\label{per1}\Lambda_k(A) \subseteq \cap\{W(A+F):F\in \BH\mbox{
has rank }< k\}\,.
\end{equation}
We show that  the inclusions in (\ref{per}) and (\ref{per1})
become inequalities if $\dim\cH$ is finite. Examples are 
given to show that these not true for infinite dimensional operators.
Nevertheless, we are able to show that
equalities also hold in (\ref{per}) and (\ref{per1})
for infinite dimensional normal operators in Section 4.
The set equalities in  (\ref{per}) and (\ref{per1}) can be viewed as 
refinements of (\ref{Linf}).  Similar set equality results are given in
Corollary \ref{closure}, which can be viewed as refinements of 
(\ref{ess}).
In Section 5, we extend the results in Sections 2 -- 4
to $\Lambda_\infty(A)$. In particular, we show in Theorem \ref{inf} and
\ref{last} that
\begin{equation} \label{1.7}
\Lambda_{\infty}(A) =\bigcap_{k\ge 1}\Lambda_k(A)
=\cap\{W(A+F): F \in \BH \hbox{ has finite rank} \},
\end{equation}
and $\Lambda_\infty(A)\ne \emptyset$ if and only if  the closure of
$\Lambda_\infty(A)$ is the essential numerical range of $A$.
Then we determine the condition under which
$\Lambda_\infty(A)$ is nonempty. The first equality
in (\ref{1.7}) gives an affirmative answer to a question
of Martinez-Avendano \cite{R}.

We close this section by listing some basic properties for the higher rank
numerical range; see \cite{Cet,Cet0,Cet1,Cet2,GLW,LPS,LS,W1}.

\begin{itemize}
\item[P1.]
For any  $a, b \in \IC$,
$\Lambda_k(aA + bI) = a\Lambda_k(A) + b$.
\item[P2.]  For any unitary  $U \in \BH$,
$\Lambda_k(U^*AU) = \Lambda_k(A)$.
\item[P3.] If $A_0$ is a compression of $A$ on a
subspace $\cH_0$ of $\cH$ such that
$\dim \cH_0 \ge k$, then
$\Lambda_k(A_0) \subseteq \Lambda_k(A)$.
\item[P4.] Suppose $\dim \cH < 2k$. The set
$\Lambda_k(A)$ has at most one element.
\item[P5.]  If $\dim \cH \ge 3k-2$ then
$\Lambda_k(A)$ is non-empty.
Otherwise, there is $B \in \BH$ such that $\Lambda_k(B) = \emptyset$.
\item[P6.]  $\Lambda_k(A)$ is always convex.
\item[P7.]
If $\dim \cH < \infty$, then $\Lambda_k(A) = \Omega_k(A)$  with
$$\Omega_k(A) = \bigcap_{\xi \in [0, 2\pi)}
\left\{\mu\in\IC: e^{i\xi}\mu +  e^{-i\xi}\bar \mu \le
\lambda_k(e^{i\xi}A + e^{-i\xi}A^*)\right\},$$ where $\lambda_k(H)$
denotes the $k$-th largest eigenvalue of the Hermitian matrix $H \in
M_n$.
\item[P8.] If $A \in M_n$ is a normal matrix with eigenvalues
$\lambda_1, \dots, \lambda_n$, then
$$\Lambda_k(A) =
\bigcap_{1 \le j_1 < \cdots < j_{n-k+1} \le n} \conv\{ \lambda_{j_1},
\dots, \lambda_{j_{n-k+1}}\}.$$
\end{itemize}

\section{The interior and closure of $\Lambda_k(A)$}
\setcounter{equation}{0}

First, we extend the definition of $\Omega_k(A)$ to infinite
dimensional operators.
For a self-adjoint operator $H$, let
$$\lambda_k(H) = \sup \{\lambda_k(X^*HX): X \hbox{ is an isometry from }
\IC^k \hbox{ to } \cH \hbox{ so that }
X^*X = I_k \}.$$
For $A \in \BH$,  let  $\Re (A)=(A+A^*)/2$ be the real part of $A$ and
$$\Omega_k(A) = \bigcap_{\xi \in [0, 2\pi)}
\left\{\mu \in \IC: \Re(e^{i\xi}\mu) \le
\lambda_k(\Re(e^{i\xi}A ))\right\}.$$

By definition, $\Omega_k(A)$ is a compact convex set. It may
be empty if $\dim \cH \le  3k-3$; see \cite[Theorem 4.7]{Cet0}.
\iffalse
In particular, it may be a
singleton, a line segment, or a convex set with non-empty interior.
\fi
In the finite dimensional case, we have $\Lambda_k(A) = \Omega_k(A)$
as noted in property (P7). Let $A = I_k \oplus \diag(1, 1/2, \dots)$. One
easily checks that $\Omega_k(A) = [0,1]$ and $\Lambda_k(A) = (0,1]$.
(See also Example \ref{eg3.5}.) Hence, property (P7) may not
hold for infinite dimensional operator $A$.

We continue to use $\cl(S)$ to denote
the closure of a set $S$ in $\IC$.
Let $\int(S)$ denote
the relative interior of $S$.
We have the following.

\begin{theorem} \label{interior}
Let $A \in \BH$ be an infinite dimensional operator, and let $k$ be
a positive integer. Then
$$\int(\Omega_k(A)) \subseteq \Lambda_k(A) \subseteq \Omega_k(A) =
\cl({\Lambda_k(A)}).$$
\end{theorem}

\it Proof. \rm
First, we establish the inclusion $\Lambda_k(A) \subseteq \Omega_k(A)$.
By \cite[Corollary 4]{LPS}, $\Lambda_k(A)$ is always non-empty.
Suppose $\mu \in \Lambda_k(A)$. Then there is an isometry $X: \IC^k
\rightarrow \cH$ such that $X^*X = I_k$ and $X^*AX  = \mu I_k$. As a
result, for any $t \in [0, 2\pi)$ we have
$$\Re (e^{it} \mu) \le \lambda_k(\Re (e^{it}A)).$$
Thus, $\mu \in \Omega_k(A)$.

\medskip
Next, we turn to the equality $\Omega_k(A) =\cl({\Lambda_k(A)})$
and the inclusion $\int(\Omega_k(A))\subseteq \Lambda_k(A)$.
By Corollary 4 in \cite{LPS}, $\Lambda_k(A)$ is non-empty.
We consider three cases.

\medskip
\bf Case 1. \rm Suppose  $\Omega_k(A)$ is a  singleton. Then
$\Lambda_k(A) = \Omega_k(A)$ because $\Lambda_k(A)$ is non-empty,
and  $\int(\Lambda_k(A)) = \int(\Omega_k(A)) = \emptyset$.

\medskip
\bf Case 2. \rm Suppose $\Omega_k(A)$ has non-empty interior in $\IC$.
Let $\mu$ be an interior point of $\Omega_k(A)$.
We may replace $A$ by $A-\mu I$ and assume that $\mu = 0$, i.e.,
$0\in \int(\Omega_k(A))$. Therefore,
\iffalse Since $0$ is an interior point of $\Omega_k(A)$, \fi
there exists $d> 0$
such that
$$\{\mu\in \IC: |\mu| \le d \} \subseteq\Omega_k(A).$$
Thus, for all $t \in [0, 2\pi)$, $\mu=de^{it}\in \Omega_k(A)$. Write $A
= H+iG$ where $H$ and $G$ are self-adjoint. Then
$$e^{-it}A + e^{it}A^* = 2(\cos t H + \sin t G).$$
Hence,
$$\lambda_k(e^{-it}A+e^{it}A^*)\ge e^{-it}\mu+e^{it}
\overline{\mu}\quad \Ra\quad  \lambda_k(\cos t H + \sin t G) \ge d\,.$$
Then, for each $t \in
[0,2\pi)$ there is $X_t: \IC^k \rightarrow \cH$ with $X_t^*X_t =
I_k$ such that $\lambda_k(\cos t X_t^*HX_t + \sin t X_t^* G X_t) >
d/2$. Furthermore, there is $\delta_t > 0$ such that for each $s\in
(t-\delta_t, t+\delta_t)$,
$$\|(\cos t X_t^*HX_t + \sin t X_t^* G X_t)
-(\cos s X_t^*HX_t + \sin s X_t^* G X_t)\| < d/4.$$
Note that $|\lambda_k(R) -\lambda_k(S)| \le \|R-S\|$ for any two Hermitian
matrices $R$ and $S$ by the Weyl's inequality; for example, see
\cite[III.2]{B}. It follows that
$$|\lambda_k(\cos t X_t^*HX_t + \sin t X_t^* G X_t)
-\lambda_k(\cos s X_t^*HX_t + \sin s X_t^* G X_t)| < d/4.$$
Consequently,
$$\lambda_k(\cos s X_t^*HX_t + \sin s X_t^* G X_t)
> \lambda_k(\cos t X_t^*HX_t + \sin t X_t^* G X_t) - d/4
> d/4.$$
Since $[0,2\pi]$ is compact, there exists a finite sequence $0 \le
t_1 < \cdots < t_m < 2\pi$ so that
$$[0,2\pi] \subseteq
\bigcup_{j=1}^m (t_j - \delta_{t_j}, t_j + \delta_{t_j}).$$ Let $A_0
= H_0 + i G_0$ be a compression of $A$ onto a subspace spanned by
the range spaces of $X_{t_1}, \dots, X_{t_m}$. Then $\lambda_k(\cos
t H_0 + \sin t G_0) > d/4$ for all $t\in [0,2\pi)$ and so $0\in
\Omega_k(A_0)$. Thus, $0\in \Lambda_k(A_0) \subseteq \Lambda_k(A)$
by Theorem 2..2 in \cite{LS}. Hence, $\int(\Omega_k(A)) \subseteq \Lambda_k(A)$
and thus $\cl(\Lambda_k(A)) = \Omega_k(A)$.

\medskip
\bf Case 3. \rm
Suppose $\Omega_k(A)$ is not a singleton and has no interior in $\IC$.
Since $\Omega_k(A)$ is a compact convex set in
$\IC$, if it is not a singleton and has no
interior in $\IC$, then it is a non-degenerate line segment.
We will show that $\Lambda_k(A)$ contains all the
relative interior points of $\Omega_k(A)$. The result will then follow.

Assume $\gamma$ is a (relative) interior point of the line
segment. By property (P1), we may assume that $[-1,1] \subseteq
\Omega_k(A) \subseteq \IR$ and $\gamma = 0$. Write $A = H+iG$ where
$H$ and $G$ are self-adjoint. Since $-1, 1 \in  \Omega_k(A)$, we
have $\lambda_k(\cos t H +\sin t G)\ge |\cos t|$ for all $t\in
[0,2\pi]$. We claim that $\lambda_k(G) = 0$. If it is not true, then
there is $\delta > 0$ such that $\lambda_k(\cos t H + \sin t G) \ge
\varepsilon > 0$ for each $t \in [\pi/2 - \delta, \pi/2 + \delta]$.
By decreasing $\varepsilon$, if necessary, we may assume that
$|\cos(\pi/2+ \delta)|=|\cos(\pi/2- \delta)|\ge\varepsilon$.
Therefore, we have $\lambda_k(\cos t H +\sin t G)\ge \varepsilon$
for all $t\in [0,\pi]$. Let $\mu =i \varepsilon$. Then, we have
$${\rm Re}(\mu e^{-it})\le\left\{ \begin{array}{rl}\varepsilon\le
\lambda_k(\cos t H +\sin t G)&\mbox{ if }t\in [0,\pi],\\&\\
0\le \lambda_k(\cos t H +\sin t G)&\mbox{ if }t\in
[\pi,2\pi].\end{array}\right.$$ Therefore, $i\varepsilon\in
\Omega_k(A)$. This contradicts that $\Omega_k(A)$ is a line segment in $\IR$.
Similarly, we can show that $\lambda_k(-G) = 0$. So, we may assume
that $G$ has operator matrix $D \oplus 0$ with
$$D = \diag(d_1, \dots, d_{p+q})$$
such that $d_1, \dots, d_p > 0$ and $d_{p+1}, \dots, d_{p+q} < 0$,
where $p < k$ and $q < k$. Let $H_0$ and $A_0$ be the compressions
of $H$ and $A$ to the kernel of $G$, respectively.

\medskip
Suppose $\lambda_k(H_0) > 0$ and $\lambda_k(-H_0) > 0$. Then $H_0$
has a compression $\tilde H_0 \in M_{2k}$ such that $\tilde H_0$ has
$k$ positive eigenvalues and $k$ negative eigenvalues. Clearly, $0
\in \Lambda_k(\tilde H_0)$ and $\tilde H_0$ is also a compression of
$A$. Then $0\in \Lambda_k(H_0) \subseteq \Lambda_k(A)$. So, we
assume that $\lambda_k(H_0) \le 0$ without loss of generality.

\medskip
Suppose the kernel of $H_0$ has dimension at least $k$. Then again
we have $0 \in \Lambda_k(H_0) = \Lambda_k(A_0) \subseteq \Lambda_k(A)$. Thus, we may
assume that the kernel of $H_0$ has dimension  $<k$. Then  $H_0$ has
operator matrix of the form
$$H_{22} \oplus H_{33}$$ so that
$H_{22} \in M_r$ with $r < 2k-1$ is positive semi-definite and
$H_{33}$ is negative definite such that the kernel of $H_{33}$ is
the zero space. Clearly,  there is a negative real number in
$\Lambda_k(H_{33}) \subseteq \Lambda_k(A_0) \subseteq \Lambda_k(A)$.
We will show that $\Lambda_k(A)$ also contains a positive real
number. By the convexity of $\Lambda_k(A)$, it will then follow that
$0 \in \Lambda_k(A)$.

\medskip
Note that $0$ is an interior point, and $H_{22}$ is finite
dimensional. We may find a small $\varepsilon > 0$ such that
$\varepsilon \in \Omega_k(A)$ and $H_0 - \varepsilon I = \hat H_{22}
\oplus \hat H_{33}$ so that $\hat H_{22}$ is positive semi-definite
and $\hat H_{33}$ is negative definite bounded above by
$-\varepsilon < 0$. Thus, $\hat H_{33}$ is invertible, and there is
an orthonormal basis of $\BH$ so that the operator matrices of $G$
and $\hat H = H - \varepsilon I $ equal
$$D\oplus 0 \qquad   \hbox{ and } \qquad
\begin{bmatrix}
\hat H_{11} & \hat H_{12} & \hat H_{13} \cr \hat H_{21} & \hat
H_{22} & 0 \cr \hat H_{31} & 0 & \hat H_{33} \end{bmatrix}$$
for $\hat H_{22} \in M_{r'}$ with $r' \le r$. For notational simplicity,
we rename $r'$ as $r$. Suppose $S \in \BH$ has operator matrix
$$\begin{bmatrix}
I_{p+q} & 0 & -\hat H_{13}\hat H_{33}^{-1}  \cr 0 & I_r & 0 \cr 0 &
0 & I \end{bmatrix}.$$ Then $SGS^*$ and $S\hat HS^*$ have operator
matrices
$$D \oplus 0 \quad \hbox{ and } \quad
\begin{bmatrix}
\hat H_{11} - \hat H_{13} \hat H_{33}^{-1} \hat H_{31} & \hat H_{12}
\cr \hat H_{21} & \hat H_{22}  \end{bmatrix} \oplus \hat H_{33}.$$
Since $0 \in \Omega_k(A-\varepsilon I)$, we see that for each $t \in
[0, 2\pi)$, we have $\lambda_k(\cos t \hat H + \sin t G) \ge 0$ and
hence $\lambda_k(\cos t S\hat HS^* + \sin t SGS^*) \ge 0$.
Consequently, if we let $\tilde H_{33}$ be the leading $k\times k$
submatrix of $\hat H_{33}$ and let $\tilde A = \tilde H + i \tilde G
\in M_{p+q+r+k}$ with
$$\tilde H = \begin{bmatrix}
\hat H_{11} - \hat H_{13} \hat H_{33}^{-1} \hat H_{31} & \hat H_{12}
\cr \hat H_{21} & \hat H_{22}  \end{bmatrix} \oplus \tilde H_{33}
\quad \hbox{ and } \quad \tilde G = D \oplus 0_{r+k},$$ then
$\lambda_k(\cos t \tilde H + \sin t \tilde G) \ge 0$ for all $t \in
[0, 2\pi)$ and hence $0 \in \Omega_k(\tilde A)$.
By Theorem 2.2 in \cite{LS}, there is a $(p+q+r+k)\times k$ matrix $X$
such that
$$X^*X = I_k \quad \hbox{ and }
\quad X^* \tilde AX = 0_k.$$
Consequently, as $\tilde A$ is a finite compression of $S(A-\varepsilon I)S^*$,
there is a partial isometry $Y: \IC^k \rightarrow \cH$ such that   $Y^*S(A-\varepsilon
I)S^*Y = 0_k$. Note that $S^*Y = ZT$ with $Z^*Z = I_k$ for some
invertible $T \in M_k$. Thus, $Z^*(A-\varepsilon I)Z = 0_k$, i.e.,
$\varepsilon \in \Lambda_k(A)$. \qed

\medskip
In the finite dimensional case, $\Lambda_k(A)$ is always closed.
If $\dim\cH$ is uncountable, then for any bounded convex set $S$ in $\IC$
one can construct a normal operator $B$ using the points in $S$ as
diagonal elements so that $\Lambda_k(A) = S$ for $A = B \otimes I$.
In the following, we  give examples of $A$
acting on a separable Hilbert space such that $\Lambda_k(A)$ has
non-empty interior with no, some or all its boundary points.
It is known that $\Lambda_k(A)$ is a singleton if $A$ is a scalar
operator, and that $\Lambda_k(A) \subseteq \IR$ if $A = A^*$.
We give examples different from these trivial cases.

\bigskip\noindent
\begin{example} \rm In the following examples,
let $B = \begin{bmatrix}0 & 2 \cr 0 & 0 \end{bmatrix}$.

\begin{itemize}
\item[(a)] Let $A = B\otimes I_{k-1} \oplus 0$. Then
$\Omega_k(A) = \Lambda_k(A) = \{0\}$.

\item[(b)]
Let $A = B\otimes I_{k-1} \oplus \diag(1,1/2, 1/3, \dots)$.
Then $\Lambda_k(A) =
(0,1]$. One can easily modify the example so that $\Lambda_k(A) =
[0,1]$ or $\Lambda_k(A) = (0,1)$.

\item[(c)]
Let $A = B\otimes I_{k-1} \oplus C$. If $C = B \oplus 0$ then
$\Lambda_k(A)$ is the closed unit disk; if $C$ is the unilateral
shift, then $\Lambda_k(A)$ is the open unit disk; if $C = \diag(-1,
i, -i, 1/2, 2/3, 3/4, 4/5, \dots)$ then $\Omega_k(A)$ is the convex
hull of $\{-1, i, -i, 1\}$, and $\Lambda_k(A)$ is the union of the
interior of $\Omega_k(A)$ and the convex hull of $\{-1,i, -i\}$.
\end{itemize}
\end{example}

\section{Low rank
perturbations of $\Lambda_k(A)$ for general  operators}
\setcounter{equation}{0}

For a positive integer $r$,
let $\cF_r$ be the set of operators in $\BH$ with rank at most $r$,
and let $\cP_r$ be the set of
rank $r$ orthogonal projections in $\BH$.

\begin{theorem}\label{inclusion}
Let $1\le r<k<\infty$. Suppose $A\in \BH$ and $F \in \cF_r$.
Then $\Lambda_{k}(A)
\subseteq \Lambda_{k-r}(A+F).$
Consequently,
$$\Lambda_{k}(A)
\subseteq \cap \{\Lambda_{k-r}(A+F): F \in \cF_r\}.$$
\end{theorem}

\it Proof. \rm
Suppose $\lambda\in\Lambda_{k}(A)$.
Let $X:\IC^k\to \cH$ be an isometry such that $X^*AX=\lambda I_k$.
Then $X^*FX$ has rank at most $r$.
There is a unitary $U \in M_k$ such that
$$U^*X^*FXU = \begin{bmatrix}
0_{k-r} & * \cr 0 & * \end{bmatrix}.$$
Let $U_1$ be obtained by taking the first $k-r$ column of $U$, and
$V = XU_1$. Then $V^*(A+F)V = \lambda I_{k-r}$ so that
$\lambda \in  \Lambda_{k-r}(A+F)$. \qed

Note that one can easily adapt the above proof to show that for
$A_1, \dots, A_m \in \BH$,
if $X^*A_jX = \lambda_j I_k$ with $X^*X = I_k$
and if $F_1, \dots, F_m \in \BH$ are such that
$$U^*X^*F_jXU = \begin{bmatrix}
0_{k-r} & * \cr
0 & * \end{bmatrix}\quad j = 1, \dots, m,$$
then $V^*A_jV = \lambda_j I_{k-r}$ for all $j = 1, \dots, m$.
So, the comment about a low rank perturbation of a quantum
channel in  Section 1 follows.

\iffalse
\it Proof. \rm Suppose $\lambda\in\Lambda_{k}(A)$.
Let $X:\IC^k\to \cH$ be an isometry such that $X^*AX=\lambda I_k$.
Suppose $F\in \BH$ has rank at most $r$.
Then there exist $u_j, v_j\in \cH$,
$1\le j\le \hat r \le r$ such that
$Fx=\sum_{j=1}^{\hat r} \la x,v_j\ra u_j$. Let $W$
be the range of $X$ and $P_W$ be the orthogonal projection of $\cH$
onto $W$. Let $V$ be the subspace of $\BH$ spanned by
$\{P_W(v_i):1\le i\le \hat r\}$. Then $\dim V\le \hat r<k$ Therefore,
$\dim (V^\perp\cap W)=m\ge k-r\ .$ Let  $\{x_1,\dots,x_m\}$ be an orthonormal
basis of $V^\perp\cap W$. Then for $1\le i,j\le m$, we have
$$\la  (A+F)x_i,x_j\ra=\la Ax_i,x_j\ra=\lambda\delta_{i\,j}.$$
Hence, $\lambda\in
\Lambda_{k-r}(A+F).$
\qed
\fi

\iffalse
\medskip
If $1 \le r < k \le n$ and $A \in M_n$, then  $\Lambda_k(A)$
can be written as the intersection of $\Lambda_{k-r}(A+F)$
for a collection of rank $r$ matrices $F$ as shown in the
following.
\fi

\medskip
If $1 \le r < k \le \dim \cH$ and $A \in \BH$, then $\Omega_k(A)$
can be written as the intersection of $\Omega_{k-r}(A+F)$
for a collection of rank $r$ operators $F$ as shown in the
following.

\begin{theorem}\label{Omega}
Suppose $A\in \BH$ and $1 \le r < k<\infty$. Let
$\cS$ be a subset of $\cF_r$ containing  the set
$\cS_0=\{2e^{i\xi}\|A\|P: P\in \cP_r \mbox{ and }\xi \in [0, 2\pi) \}.$
Then
$$\Omega_k(A) = \cap\{\Omega_{k-r}(A+F): F \in \cS\}.$$
\end{theorem}

\it Proof. \rm 
The inclusion ($\subseteq$) follows from Theorem \ref{inclusion}
and the fact that $\Omega_k(A) = \cl (\Lambda_k(A))$ by Theorem \ref{interior}.

Suppose $\lambda\notin \Omega_k(A)$. 
Then there exists $t\in \IR$ such that
$\lambda_k(\Re\(e^{it}A\)) < \Re\(e^{it}\lambda\)$.
Let $e^{it}A=H+iG$ be with $H = H^*$ and $G=G^*$.
Then $H$ has an operator matrix
$\diag(\lambda_1, \dots, \lambda_m) \oplus H_2$
with $m \le k-1$
such that $\sup \sigma(H_2) < \Re(e^{it}\lambda)$.
Let $F= -2e^{-it}\|A\|(I_m \oplus I_{k-m} \oplus 0)\in \BH$.
Then $\lambda_{k-r} (\Re(e^{it}(A+F)) < \Re\(e^{it}\lambda\)$.
Hence,
$\lambda \notin \Omega_{k-r}(A+F)$.
\qed

\iffalse
\it Proof. \rm
The inclusion ($\subseteq$) follows from Theorem \ref{inclusion}
and the fact that $\Omega_k(A) = \cl (\Lambda_k(A))$ by Theorem \ref{interior}.

Suppose $\lambda\notin \Omega_k(A)$. We may assume that $\lambda=0$.
Then there exists $t\in \IR$ such that
$\lambda_k(\Re\(e^{it}A\))<0$. Hence, $A\ne 0$.
Let $e^{it}A=H+iG$ be with $H = H^*$ and $G=G^*$.
Then $H$ has an operator matrix
$\diag(\lambda_1, \dots, \lambda_m) \oplus H_2$
with $m \le k-1$
such that $\sup \sigma(H_2) < 0$.
Let $F= -2e^{-it}\|A\|(I_m \oplus I_{k-m} \oplus 0)\in \BH$.
Then  $\Re(e^{it}(A+F))\le -\|A\|I$.
Hence,
$\lambda \notin \Omega_{k-r}(A+F)$.
\qed
\fi

Note that for the set $\cS$ in the above theorem, we can take the whole   
$\cF_r$ or the much smaller subset $\cS_0$.
\iffalse
With Theorem \ref{interior}, we have the following corollary.

\begin{corollary} \label{closure}
Under the same setting as in Theorem \ref{Omega}.
Each of the following sets is equal to $\cl( \Lambda_k (A) )$.
\begin{itemize}
\item[\rm (a)] $\cap \{\cl( \Lambda_{k-r}(A+F)): F \in \cS\}$.
\item[\rm (b)] $\cap \{\cl( \Lambda_{k-1}(A+ 2e^{i\xi}\|A\|P) ):
\xi \in [0, 2\pi),\ P\in \cP_1\}$.
\item[\rm (c)] $\cap \{\cl(W(A+2e^{i\xi}\|A\|P)):
\xi \in [0, 2\pi),\ P\in \cP_{k-1} \}$.
\end{itemize}

\end{corollary}
\fi
We have the following corollary.

\begin{corollary} \label{closure}
Under the same setting as in Theorem \ref{Omega}.
Each of the following sets is equal to $\Omega_k (A) $.
\begin{itemize}
%\item[\rm (a)] $\cap \{\Omega_{k-r}(A+F): F \in \cS\}$.
\item[\rm (a)] $\cap \{\Omega_{k-1}(A+ 2e^{i\xi}\|A\|P) :
\xi \in [0, 2\pi),\ P\in \cP_1\}$.
\item[\rm (b)] $\cap \{\Omega_1(A+2e^{i\xi}\|A\|P):
\xi \in [0, 2\pi),\ P\in \cP_{k-1} \}$.
\end{itemize}
\end{corollary}
With Theorem \ref{interior}, the above result also holds if we replace $\Omega_m(B)$ by $\cl\(\Lambda_m(B)\)$.
Using the fact that $\Lambda_k(A) = \Omega_k(A)$
when $A\in M_n$, we have the following result.

\begin{theorem}\label{finite}
Suppose $A\in M_n$ and $1\le r<k \le n$.  Let
$\cS$ be a subset of $\cF_r$ containing the set
$\cS_0 = \{2e^{i\xi}\|A\|P: P\in \cP_r \mbox{ and }\xi \in [0, 2\pi) \}$.
Then
\begin{itemize}
\item[\rm (a)] $\Lambda_k(A) = \cap \{\Lambda_{k-r}(A+F): F \in \cS\}$.
\item[\rm (b)] $\Lambda_k(A) = \cap \{\Lambda_{k-1}(A+ 2e^{i\xi}\|A\|P):
\xi \in [0, 2\pi), P\in\cP_1\}$.
\item[(c)] $\Lambda_k(A) = \cap \{W(A+2e^{i\xi}\|A\|P):
\xi \in [0, 2\pi),\ P\in \cP_{k-1} \}.$
\end{itemize}
\end{theorem}

\iffalse
\begin{theorem}\label{Omega} Suppose $A \in M_n$.
For $1\le r < k \le n$, let $\cS$ be a subset of $\cF_r$
containing the set
\ $\{2e^{i\xi} \|A\| P: \xi \in [0, 2\pi), \  P\in \cP_r \}.$
Then
$$\Lambda_k(A) = \cap \{\Lambda_{k-r}(A+F): F \in \cS\}.$$
\end{theorem}

\it Proof. \rm  The inclusion
``$\subseteq$'' follows from Theorem \ref{inclusion}.

Suppose $\lambda\not\in\Lambda_{k}(A)=\Omega_k(A)$. Then there
exists $\xi\in \IR$ such that $\Re(e^{i\xi}\lambda ) >
\lambda_k\(\Re(e^{i\xi}A )\)$. Let
$\Re(e^{i\xi}A )=U^*\diag\(a_1,\dots,a_n\)U$, where $a_1\ge \dots
\ge a_n$. Let $P= U^*\(I_r\oplus 0_{n-r}\)U$.
Then $P\in \cP_r$ and
$$ \Re(e^{i\xi}(A-2\|A\|P))
= U^*\diag \(a_1-2\|A\|,\dots,a_r-2\|A\|, a_{r+1},\dots, a_n\)U.$$
So,
$\Re(e^{i\xi}\lambda)>a_k=\lambda_{k-r}(\Re(e^{i\xi}(A-2\|A\|P)))$
and hence
$$\lambda\not \in \Lambda_{k-r}(A- e^{-i\xi}2\|A\| P).$$
\vskip -.3in\qed

\

\medskip

\begin{corollary} Let $A\in M_n$ and $k\ge 1$. We have
\begin{itemize}
\item[(a)] $\Lambda_{k}(A)=\cap \{\Lambda_{k-1}(A+ 2e^{i\xi}\|A\|xx^*):
\xi \in [0, 2\pi), \ x\in\IC^n,\ x^*x = 1\}$.
\item[(b)] $\Lambda_k(A) =
\cap \{W(A+2e^{i\xi}\|A\|P):
\xi \in [0, 2\pi),\ P\in \cP_{k-1} \}.$
\end{itemize}
\end{corollary}
\fi

\medskip
The following example shows that Theorem \ref{finite}
does not hold for infinite dimensional operators.

\begin{example}\label{eg3.5} Let $A = A_1\oplus A_2$,
where
$$A_1 = \[\begin{array}{cc}0 & i \\ i & 2\end{array}\]
\quad \mbox{ and }  \quad
A_2 =  {\rm diag}(b_2, \bar b_2, b_3, \bar b_3, \cdots) \oplus
\, {\rm diag}(b_2, \bar b_2, b_3, \bar b_3, \cdots)$$
with $b_m = -1+ e^{i\pi/m}$  for $m =  2,
\dots$.  Then $0 \in \cl({\Lambda_2(A)})$
and  $0 \notin \Lambda_2(A)$,
but $0 \in \cap\{ W(A+F): F \hbox{ is rank one} \} $.

\medskip
\it Verification. \rm
Note that every $\mu \in \Lambda_1(A_{2})$
is an element of $\Lambda_2(A_{2})$, and hence
$\Lambda_1(A_{2}) = \Lambda_2(A_{2}).$
Clearly, $0 \in \cl(\Lambda_1(A)) = \cl(\Lambda_2(A))$.

Next, we show that $0\notin \Lambda_2(A)$.
Suppose $0 \in \Lambda_2(A)$. Then $0 \in \Lambda_2(H)$ for $H =
(A+A^*)/2$. Let $U$ be unitary such that $U^*AU =
\[\begin{array}{cc}0_2 & * \\ * & * \end{array}\]$.
Then $U^*HU$ has the same form. Since $H$ has spectrum $\{2,0\} \cup
\{-1+\cos \pi/m: m =  2, \dots\},$ we may assume that $U$ has the
form $[1] \oplus U_1$ such that the $(1,1)$ entry of $U_1$ is
nonzero. But then $U^*GU$ will have non-zero $(1,2)$ entry for $G =
(A-A^*)/(2i)$. This contradicts the fact that $U^*AU$ has zero
$(1,2)$ entry. So, we see that $0 \notin \Lambda_2(A)$.

\medskip
Now, suppose $F = \[\begin{array}{cc}F_{11}& F_{12}\\
F_{21} & F_{22} \end{array}\]$
is a rank one operator with $F_{11}\in M_2$.
Let $x \in \IC^2$ be a nonzero vector such that $F_{11}x = 0$.
If $x$ is  a multiple of $e_1$, then the $(1,1)$ entry
of $A_{1} + F_{11}$ equals $0$ and we have $0 \in W(A_{1}+F_{11}) \subseteq
W(A+F)$. If $x$ is not a multiple of $e_1$,
then $\mu_0 = x^*(A_1 + F_{11}) x = x^* A_1 x\in W(A_1)$ has positive real part
and $\mu_0\in W(A_{1}+F_{11}) \subseteq W(A+F)$.
Since $F_{22}$ has rank at most one, by Theorem \ref{inclusion} we have
$$W(A_{2}) = \Lambda_2(A_{2})
\subseteq \Lambda_1(A_{2}+F_{22}) = W(A_{2}+F_{22}) \subseteq W(A + F).$$
So there exist
$\mu_1, \mu_2 \in W(A_{2}) \subseteq W(A+F)$
on the different sides of the line passing through $\mu_0$ and the origin.
It follows that $0 \in \conv\{\mu_0, \mu_1, \mu_2\} \subseteq W(A+F)$
by the convexity of $W(A+F)$.
Consequently, we have
$$0 \in \cap\{ W(A+F): F \hbox{ has rank one} \}.$$
\vskip -.3in \qed
\end{example}

\section{Low rank
perturbations of $\Lambda_k(A)$ for infinite dimensional normal operators}
\setcounter{equation}{0}

\medskip
In the following, we prove that Theorem \ref{finite} is valid for
(infinite dimensional) normal operators.
We first establish some auxiliary results showing
that one can refine the spectral decomposition of
a normal operator using the geometrical information of
its numerical range.

\iffalse
%For $A\in \BH$ and $\theta\in \IR$, let
$A(\theta)=e^{i\theta}A$, $H(\theta)=(A(\theta)+A(\theta)^*)/2$,
%$G(\theta)=(A(\theta)-A(\theta)^*)/(2i)$.

Suppose
$A=H+iG$ with $H$, $G$ self-adjoint. Define $\mu_k(A)=\lambda_k(G)$ and
$\mu_k(A,\theta)=\mu_k\(e^{-i\theta}A\)$, for $\theta\in\IR$.
%and  $\cL=\{z\in \IC:\Re(z)<0\}$ the left half plane.
\fi

Let $\cP=\{z\in \IC:\Im(z)>0\}$ be the open upper half plane of $\IC$.
For $A \in \BH$ and $k \le \dim \cH$, let
$$\mu_k(A, t) = \lambda_k((e^{-it} A - e^{it}A^*)/(2i)).$$
Notice also that
$$\Omega_k(A) = \bigcap_{t \in [0,2\pi)}
\left\{\mu\in \IC: \Im (e^{-it} \mu) \le \mu_k(A,t) \right\}.$$

\begin{lemma}\label{decomp} Suppose $A\in \BH$ is normal.
If $\mu_m(A,t)\le 0$ for some $m\ge 1$ and $t\in \IR$.
Then $A$ has a decomposition
$A_1\oplus A_2\oplus\hat A$ such that
$\dim A_1  < m$,
$$W(A_1)\subseteq e^{it}\cP, \quad
W(A_2)\subseteq -e^{it}\cP \quad \hbox{ and } \quad
W(\hat A)\subseteq e^{it}\IR.$$
Furthermore, if $\lambda_\ell( e^{-it}\hat A+e^{it} \hat A^*)/2 \le 0$
for some $\ell \ge 1$, then $\hat A$  has a decomposition $A_3 \oplus  A_4\oplus 0$
such that dim $A_3<\ell$,
$$W(A_3)\subseteq e^{it}(0,\infty),
\quad \hbox{ and } \quad
W(A_4) \subseteq e^{it}(-\infty,0).$$
Note that each of the summands $A_1, A_2, \hat A, A_3, A_4$
may be vacuous.
\end{lemma}

\it Proof. \rm Without loss of generality, we may assume that $t=0$.
Let  $A=H+iG$, where $H,G$ are self-adjoint. Then
$G = G_1 \oplus G_2\oplus 0$ such that $G_1$
is positive definite with dimension $p < m$
and $G_2$ is negative definite.
Let
$$H =
\begin{bmatrix}H_{11} & H_{12} & H_{13} \cr
H_{12}^* & H_{22}&H_{23}\cr H_{13}^* & H_{23}^* & H_{33} \end{bmatrix}$$
such that $H_{12} = [D \, | \, 0]$, where
$D = \diag(d_1, \dots, d_p)$ with $d_1 \ge \cdots \ge d_p \ge 0$.
Since $ GH=HG$,
it follows that $G_1 [D \, | \, 0] = [D \, | \, 0]G_2$.
Since $G_1$ is positive definite and $G_2$ is negative definite,
the $(1,1)$ entry on the left side is nonnegative,
and the $(1,1)$ entry on the right side is nonpositive.
Thus, $d_1 = 0$ and hence $H_{12} = 0$.
Since $G_1H_{13}=0$  and $G_2H_{23}=0$, we have $H_{13}=0$  and $H_{23}=0$.
So, $H=H_{11}\oplus H_{22}\oplus H_{33}$ and $A$  has
asserted properties, with $A_1=H_{11}+iG_1$, $A_2=H_{22}+iG_2$,
and $\hat A=H_{33}$.

If $\lambda_\ell( e^{-it}\hat A+e^{it} \hat A^*)/2 \le 0$ for some $\ell$,
then we can apply the above result to $\hat A$ and get the
desired decomposition for $\hat A$.
\qed

The following result \cite[Lemma 2 and Corollary]{AS} will be needed in later discussion.

\begin{lemma}\label{LAS} Suppose dim $\cH$ is infinite. Let $T\in\BH$. Then the following are equivalent.
\begin{itemize}
\item[\rm (a)] $\lambda\in W_e(T)$.
\item[\rm (b)] There is an orthonormal set $\{e_n\}$ such that $\la Te_n,e_n\ra\to \lambda$.
\item[\rm (c)] There is a decomposition of $\cH$ as $\cH_1\oplus \cH_2$ and a sequence $\{\lambda_i\}$ in $\IC$, such that $\lambda_i\to \lambda$ and
$$T=\[\begin{array}{ccc|c} \lambda_1& & 0&\\
&\lambda_2&&*\\
0&&\ddots&\\
\hline
&*&&*\end{array}\]$$\end{itemize}
Furthermore, if $\alpha$, $\beta\in  W_e(T)$, then there exist two sequences $\{\alpha_i\}$ and $\{\beta_i\}$ in $\IC$, such that $\alpha_i\to \alpha$, $\beta_i\to \beta$ and a decomposition of $\cH$ as $\cH_1\oplus \cH_2$ such that
$$T=\[\begin{array}{ccccc|c} \alpha_1& & &&&\\
&\beta_1&& 0&& \\
&&\alpha_2&&&*\\
&0&&\beta_2&&  \\
&&&&\ddots&\\
\hline
&&&*&&*\end{array}\]$$
In both cases, we may take $\cH_2$  to be infinite dimensional.
\end{lemma}

\iffalse
Let $\cR=\{z\in\IC:\Re z>0\}$. Suppose $B\in M_n$. For each $\theta\in [0,
2\pi)$, We can apply the above Lemma  and get a decomposition of
$B=B_1(\theta)\oplus B_2(\theta)\oplus \hat B(\theta)$ such that
$$W\(B_1(\theta)\)\subseteq e^{i\theta}\cR,\ W\(B_2(\theta)\)\subseteq
e^{i\theta}\cR,\ W\(\hat B_3(\theta)\)\subseteq e^{i\theta}\IR $$

Then by the result in \cite{LS}, we have $0\in\Lambda_k(B)$ if and only if dim
$\(B_1(\theta)\oplus\hat B(\theta)\ge k$ for all $\theta\in [0,2\pi)$.
\fi

\begin{lemma}\label{LR}
Suppose $T\in \BH$ is a normal operator such that for some
$\pi\le s_1<s_2\le 2\pi$,
$$\sigma(T) \subseteq \{\rho e^{it} \in \IC: \rho > 0,\  t \in [s_1,s_2]\}.$$
Let $k$ be a positive integer and $s_3 \in (s_1,s_2)$,
$$
\cL = \{\rho e^{i t} \in \IC: \rho > 0,\ t \in (s_1,s_3)\}
\quad\hbox{and}\quad
\cR = \{\rho e^{i t} \in \IC: \rho > 0,\ t \in (s_3,s_2)\}.$$
We have
\begin{itemize}
\item[\rm (a)] If $\cL \cap \sigma(T)$ is infinite
or contains an eigenvalue of T with infinite multiplicity,
then  $T$ has a  compression  $T_1\in M_k$ such that
$W\(T_1\) \subseteq \cL$.
\item[\rm (b)] If $\cR \cap \sigma(T)$ is infinite
or contains an eigenvalue of T with infinite multiplicity,
then  $T$ has a compression  $T_2\in M_k$ such that
$W\(T_2\) \subseteq \cR$.
\end{itemize}
If both   hypotheses in (a) and (b) hold,
then $T$ has a compression of the form
$T_1 \oplus T_2$ such that  $\dim T_1 = \dim T_2 = k$
and
\begin{equation*} \label{eq1}
W(T_1) \subseteq \cL \quad \mbox{ and } \quad W(T_2) \subseteq \cR.
\end{equation*}
\end{lemma}

\iffalse
%Furthermore, let $$\cL = \{\rho e^{i t} \in \IC: \rho > 0,
%t \in (s_1,s_1+\delta)\}\,.$$
%If $\cL \cap\sigma(T)$ is also infinite, then $T$ has
%a compression $T_1\oplus T_2$ such that
%$T_1, T_2 \in M_k$ are normal satisfying (\ref{eq1}) and
%\begin{equation} \label{eq2}
%W(T_2) \subseteq \cL.
%\end{equation}
\fi

\it Proof. \rm We will prove the last assertion. The proof of (a) and (b) are similar.
Suppose both $\cR \cap \sigma(T)$ and $\cL \cap \sigma(T)$
contain only isolated points of  $\sigma(T)$.
Then we can construct $T_1$ (respectively, $T_2$)
from any  $k$ (counting multiplicity) eigenvalues of $T$ in $\cL$
(respectively, in $\cR$)
and the corresponding eigenvectors.

Suppose one of the sets $\cL \cap \sigma(T)$ or $\cR \cap \sigma(T)$, say,
$\cL \cap \sigma(T)$, contains only isolated points of
$\sigma(T)$, and the other set contains an  accumulation point of $\sigma(T)$.
Then we can construct $T_1$ from any  $k$ eigenvalues of $T$ in $\cL$ and the
corresponding eigenvectors. Let $\cH_1$ be the $k$-dimensional subspace
spanned by the $k$-eigenvectors. Then with respect to the decomposition
$\cH=\cH_1\oplus \cH_1^{\perp}$, $T=T_1\oplus S$ for some normal $S$.
Since
$\cR$ contains an accumulation point of  $\sigma(S)$ and $S$ is normal,  by
Lemma \ref{LAS}, $S$ has a $k$-dimensional 
compression $T_2$ with  $W(T_2) \subseteq \cR$.

Finally, suppose both $\cL$ and $\cR$ contain an accumulation point
of $\sigma(T)$. Then
the result follows from the last statement in Lemma \ref{LAS}.
\qed

\iffalse
\begin{theorem} \label{notin}
Suppose $A \in B(H)$ is normal,
and $\lambda \in \Omega_k(A) \setminus \Lambda_k(A)$.
Then  $A$ can be decomposed into $\tilde A_1 \oplus \tilde A_2$ such that
$\tilde A_1$ has dimension at most $k-1$,
$W(\tilde A_1) \subseteq \lambda + S$ and $W(\tilde A_2)
\subseteq \IC\setminus(\lambda +S)$,
where $S =e^{i\theta}\( \{z: {\rm Im}(z) > 0\}\cup \tilde L\)$
with $\tilde L = (-\infty,0]$ or $[0,\infty)$ for some $\theta \in [0, 2\pi)$.
Consequently, for any $1 \le r < k$
there exists an orthogonal projection $P\in \cP_r$
such that
$$\lambda \notin \Lambda_{k-r}(A - 2e^{i\theta}\|A\|P).$$
\end{theorem}
\fi

\begin{theorem} \label{notin}
Suppose $A \in B(H)$ is normal.
Then $\lambda \notin \Lambda_k(A)$
if and only if
$A$ can be decomposed into $\tilde A_1 \oplus \tilde A_2$ such that
$\tilde A_1$ has dimension at most $k-1$,
$W(\tilde A_1) \subseteq \lambda + S$ and $W(\tilde A_2)
\subseteq \IC\setminus(\lambda +S)$,
where $S =e^{it}\( \cP \cup \tilde L\)$
with $\cP = \{z\in \IC:\Im(z) >0\}$ and
$\tilde L = (-\infty,0]$ or $[0,\infty)$ for some $t \in \IR$.
%Consequently, for any $1 \le r < k$
%there exists an orthogonal projection $P\in \cP_r$
%such that
%$$\lambda \notin \Lambda_{k-r}(A - 2e^{i\theta}\|A\|P).$$
\end{theorem}

\it Proof. \rm
Suppose $A$ has the decomposition as stated
with $\dim \tilde A_1 = m \le k-1$. Take $F\in M_m$
such that $W(\tilde A_1 + F) \subseteq \IC \setminus (\lambda + S)$.
By Theorem \ref{inclusion},
$$\Lambda_k(A) \subseteq  W(A + (F \oplus 0) )
= W( (\tilde A_1 + F) \oplus \tilde A_2 )
%\subseteq \conv W( \tilde A_1 + F) \cup W(\tilde A_2 )
\subseteq \IC \setminus (\lambda + S).$$
Hence, $\lambda \notin \Lambda_k(A)$.

Conversely, suppose $\lambda \notin \Lambda_k(A)$.
Without loss of generality, we may assume that $\lambda=0$.

\bf Case 1. \rm Suppose $\lambda = 0 \notin \Omega_k(A)$, then
$\mu_k(A,t) < 0$ for some $t\in [0,2\pi)$.
By Lemma \ref{decomp}, $A = A_1 \oplus A_2 \oplus \hat A$ with
$\dim A_1 < k$, $W(A_1)\subseteq e^{it}\cP$, $W(A_2)\subseteq -e^{it}\cP$
and $W(\hat A)\subseteq e^{it}\IR$.
Furthermore, as $\mu_k(A,t) < 0$, we must have $\dim A_1 + \dim \hat A < k$.
Then $\hat A = A_3 \oplus A_4$ so that $W(A_3) \subseteq e^{it}[0,\infty)$
and $W(A_4) \subseteq - e^{it} (0,\infty)$.
Then the result follows with $\tilde A_1 = A_1 \oplus A_3$ and $\tilde A_2 = A_2 \oplus A_4$.

\medskip
\bf Case 2. \rm Suppose $\lambda = 0 \in \Omega_k(A)$
and such decomposition mentioned in the theorem does not exist.
Suppose $\ker A$, the kernel of $A$, has dimension $p < k$.
We may assume that $p = 0$. Otherwise, replace $A$ by the
compression of $A$ on $\(\ker A\)^\perp$ and replace $k$ by $k-p$.
We are going to derive a contradiction by showing
that $A$ has a finite dimensional compression
$B$ such that $0\in \Omega_k(B)=\Lambda_k(B) \subseteq \Lambda_k(A)$.

To construct the matrix $B$, we first show that
there exist $s_1 \le 0 \le \pi \le s_2$ with $s_2 - s_1 \le 2\pi$
such that $A=A_1\oplus A_2\oplus A_3 \oplus A_4$,
where
\begin{multline}\label{cond}
\dim   A_1 < \infty,\quad
W(A_1)\subseteq \{\rho e^{it}: \rho>0,\ t \in (s_1,s_2)\},\cr
W(A_2)\subseteq \{\rho e^{it}: \rho>0,\ t \in (s_2,s_1+2\pi)\},\cr
W(A_3)\subseteq e^{is_1} (0,\infty),\quad\hbox{and}\quad
W(A_4)\subseteq e^{is_2} (0,\infty).
\end{multline}
Then we show that $A_2\oplus A_3 \oplus A_4$ has a finite dimensional
compression $B_2 \oplus B_3 \oplus B_4$ such that $B= A_1\oplus B_2 \oplus B_3 \oplus B_4$
has $0\in \Omega_k(B)$.

\medskip
Since $0\in\Omega_k(A)$, we have $\mu_k(A,t)\ge 0$
for all $t \in [0, 2\pi)$.
If $\mu_k(A,t) > 0$ for all $t \in [0, 2\pi)$,
then $0$ lies in the interior of $\Omega_k(A)$. Hence, $0 \in \Lambda_k(A)$.
So, we may assume that there is $t_0 \in [0, 2\pi)$ such that
$\mu_k(A,t_0) = 0$. We may further assume that $t_0 = 0$.

As $\mu_k(A,t_0) = 0$, $A$ has at most $k-1$ eigenvalues in the open upper half plane.
Suppose these eigenvalues have arguments $0<t_1\le t_2\le \cdots\le t_p<\pi$, $p<k$.
Take $t_{p+1}=\pi $. Let
$g \in \{1,\dots, p+1\}$ be the smallest integer
such that $\mu_m(A,t_g-\pi) = 0$ for some $m$
and  $h \in \{0,1, \dots, p\}$ be the largest integer
satisfying  $\mu_m(A,t_h) = 0$ for some $m$.
Let $s_1=t_g-\pi$ and
$s_2=t_h+\pi$. We are going to find $A_j$ for $j = 1,2,3,4$ 
satisfying (\ref{cond}).

By Lemma \ref{decomp} with $t = t_h$,
we have
$A=\hat A_1\oplus \hat A_2\oplus \hat A$ such that
$$\dim \hat A_1 < \infty,\quad
W( \hat A_1 )\subseteq e^{it_h}\cP,\quad
W( \hat A_2 )\subseteq - e^{it_h}\cP,\quad\hbox{and}\quad
W ( \hat A )\subseteq e^{it_h}\IR.$$

Let $\hat H=\(e^{-it_h} \hat A+ e^{-it_h}\hat A^*\)/2 $.
If both $\lambda_k (\hat H )$ and $\lambda_k( -\hat H )$
are nonnegative, then we have $0\in \Lambda_k(\hat H )$,
which implies that $0\in \Lambda_k(A)$, a contradiction.
So, we have  either $\lambda_k(\hat H )$ or
$\lambda_k(-\hat H )$ is negative.
By Lemma \ref{decomp} and the assumption that ker $A=0$,
we have $\hat A = \hat A_3 \oplus \hat A_4$ with
$$W(\hat A_3) \subseteq e^{it_h}(0,\infty)\quad\hbox{and}\quad
W(\hat A_4) \subseteq -e^{it_h}(0,\infty).$$

If $t_h = t_g - \pi$, we take $A_j = \hat A_j$ for $j = 1,2,3,4$.
Then $A_1,A_2,A_3,A_4$ satisfy (\ref{cond}) with $s_1 = t_g - \pi = 0$ and
$s_2 = t_h +\pi= \pi$.

\medskip
Suppose $t_h > t_g - \pi$. Then $\dim \hat A_3$ is finite.
We further apply Lemma \ref{decomp} to $\hat A_2$ with $t = t_g - \pi$, we have
$\hat A_2= A_1'\oplus A_2'\oplus A_3'$, with
\begin{multline*}
\dim A_1'< \infty,\quad
W(A_1')\subseteq \{\rho e^{it}:\rho>0,\ t \in (t_g-\pi,t_h)\}, \cr
W(A_2')\subseteq \{\rho e^{it}:\rho>0,\ t \in (t_h - \pi, t_g  - \pi)\},\hbox{ and }
W(A_3')\subseteq e^{i (t_g-\pi)}(0,\infty).
\end{multline*}
Note that $A_4'$ is vacuous because $t_h\le t_g\le t_h + \pi$.
Then $A_1=\hat A_1\oplus \hat A_3\oplus A_1'$,
$A_2 = A_2'$, $A_3 = A_3'$, and $A_4 = \hat A_4$
will satisfy (\ref{cond}) with $s_1 = t_g -\pi$ and $s_2= t_h+\pi$.

\medskip
Now we choose a finite dimensional compression $B_2\oplus B_3\oplus B_4$
of $A_2\oplus A_3\oplus A_4$ and
show that $\mu_k(A_1 \oplus B_2\oplus B_3\oplus B_4, t) \ge 0$
for all $t\in [0,2\pi)$.
Observe that
\begin{eqnarray}\label{mu1}
\mu_k(A_1,t) \ge 0
\quad\hbox{for all}\quad t_g - \pi < t < t_h.
\end{eqnarray}

Let $B_3$ be a $k$-dimensional compression of $A_3$,
if $\dim A_3$ is infinite and $B_3=A_3$, otherwise.
We claim that
\begin{eqnarray}\label{mu2}
\mu_k(A_1 \oplus B_3,t) \ge 0\quad\hbox{for all}\quad
t_{g-1} -\pi \le t \le t_g - \pi.
\end{eqnarray}
The claim is clear if $\dim A_3$ is infinite.
Suppose $\dim A_3$ is finite and
$\mu_k(A_1 \oplus A_3,t) < 0$
for some $t \in [t_{g-1} -\pi, t_g - \pi]$.
Since $\dim(A_1\oplus A_3)$ is finite and
$W(A_1 \oplus A_3) \subseteq
\{\rho e^{it}:\rho>0,\ t \in [t_g - \pi, t_h+\pi) \}$,
$A_1 \oplus A_3$ has a decomposition $A_1'' \oplus A_3''$,
with
$$\dim A_1'' < k,\quad W(A_1'') \subseteq S,\quad
\hbox{and}\quad W(A_3'') \subseteq \IC\setminus S,$$
where $S = e^{i(t_g - \pi)}(\cP \cup [0,\infty))
= \{\rho e^{it}:\rho>0,\ t\in [t_g - \pi, t_g)\} $.
Notice also that $W(A_2 \oplus A_4) \subseteq \IC\setminus S$.
Then if we take $\tilde A_1 = A_1''$ and $\tilde A_2 = A_2 \oplus A_3'' \oplus A_4$,
we have $W(\tilde A_1) \subseteq S$ and
$W(\tilde A_2) \subseteq \IC\setminus S$,
which contradicts our assumption that such decomposition does not exist.

Next, let $B_4$ be a $k$-dimensional compression of $A_4$,
if $\dim A_4$ is infinite and $B_4=A_4$, otherwise.
By a similar argument as in the previous paragraph, we can show that
\begin{eqnarray}\label{mu3}
\mu_k(A_1 \oplus B_4,t) \ge 0
\quad\hbox{for all}\quad
t_h \le t \le t_{h+1}.
\end{eqnarray}

In the following, we will choose a finite dimension compression of $B_2$ of $A_2$ so that
\begin{eqnarray}\label{mu4}
\mu_k(B_2\oplus B_3 \oplus B_4,t) \ge 0
\quad\hbox{for all}\quad
t_{h+1} \le t \le t_{g-1} + \pi.
\end{eqnarray}

Suppose $\dim A_2$ is finite.
Then by the definition of $t_h$ and $t_g$, both $\dim A_3$ and $\dim A_4$ are infinite.
Then $\mu_k(B_3\oplus B_4,t)\ge 0$ for all $t\in [t_{h+1},t_{g-1} + \pi]$
and so (\ref{mu4}) holds with vacuous $B_2$.

\medskip
Now suppose $\dim A_2$ is infinite. We consider the following three cases.

\noindent\bf Case 1. $t_g = t_h$. \rm
In this case, the summand $A_2'\oplus A_3'$ is vacuous
and so as $A_2\oplus A_3$. Also $\dim A_4$ is infinite.
Then (\ref{mu4}) holds with vacuous $B_2$ and $B_3$.

\noindent\bf Case 2. $t_g = t_{h+1}$. \rm
Let $B_2$ be a $k$-dimensional compression of $A_2$.
Then $\mu_k(B_2,t) \ge 0$ for all $t\in[t_{h+1},t_{g-1} - \pi]$
and so (\ref{mu4}) holds.

\noindent\bf Case 3. $t_g > t_{h+1}$. \rm
Because of the choice $t_g$ and $t_h$, both
$\cL\cap \sigma(A_2)$ and $\cR\cap \sigma(A_2)$
are infinite or contains an eigenvalue of $A_2$ with
infinite multiplicity, where
$$\cL = \{\rho e^{it}:\rho>0,\ t \in (t_h+\pi,t_{h+1}+\pi)\}\hbox{ and }
\cR = \{\rho e^{it}:\rho>0,\ t \in (t_{g-1}+\pi,t_g+\pi)\}.$$
By Lemma \ref{LR}, we can get finite dimensional compressions $T_1$ and
$T_2$ of $A_2$ such that $\dim(T_1) = \dim(T_2) = k$,
$W\(T_1\)\subseteq \cL$ and $W\(T_2\)\subseteq \cR$.
Then $\mu_k(T_1,t) \ge 0$ for all $t \in [t_{h+1}, t_h+\pi]$
and $\mu_k(T_2,t) \ge 0$ for all $t \in [t_g, t_{g-1}+\pi]$.
Thus, $B_2=T_1\oplus T_2$ will satisfy (\ref{mu4}).

\medskip
Now let $B = A_1 \oplus B_2\oplus B_3\oplus B_4$. By
(\ref{mu1}), (\ref{mu2}), (\ref{mu3}), and (\ref{mu4}),
we conclude that $\mu_k(B,t) \ge 0$ for all $t\in [0,2\pi)$
and hence $0\in \Omega_k(B) = \Lambda_k(B)$.
\qed

\begin{theorem}\label{normal}
Suppose $A\in \BH$ is normal and $1 \le r < k$. Let
$\cS$ be a subset of $\cF_r$ containing  the set
$\cS_0=\{2e^{i\xi}\|A\|P: P\in \cP_r \mbox{ and }\xi \in [0, 2\pi) \}.$
Then
$$\Lambda_k(A) = \cap\{\Lambda_{k-r}(A+F): F \in \cS\}.$$
\end{theorem}

\it \it Proof. \rm The inclusion ($\subseteq$) follows from Theorem
\ref{inclusion}.
Suppose $\lambda\notin\Lambda_{k}(A)$.
By Theorem \ref{notin}, $A$ has a decomposition $A_1 \oplus A_2$
with $A_1\in M_m$, $W(A_1) \subseteq \lambda + S$ and $W(A_2) \subseteq \IC\setminus (\lambda +S)$,
where $m\le k-1$ and $S$ is defined as in Theorem \ref{notin}.
Let $F = -2ie^{it}\|A\|(I_r \oplus 0)\in \cS_0$.
Then $A+F$ has less than $k-r$ eigenvalues in $\lambda + \cS$.
Thus, $A+F$ has a decomposition $B_1 \oplus B_2$ with $\dim B_1 < k-r$
such that $W(B_1) \subseteq \lambda + S$ and $W(B_2) \subseteq \IC\setminus (\lambda +S)$.
By Theorem \ref{notin}, $\lambda \notin \Lambda_{k-r}(A+F)$.
\qed

\iffalse
\it \it Proof. \rm The inclusion ($\subseteq$) follows from Theorem
\ref{inclusion}.
Suppose $\lambda\notin\Lambda_{k}(A)$. If $\lambda\notin \Omega_k(A)$, then
by Theorem \ref{Omega}, there is $F \in \cS$ such that $\lambda\notin \Omega_{k-r}(A+ F)$.
Hence, $\lambda \notin \Lambda_{k-r}(A+F)$.
If $\lambda\in \Omega_k(A)$, then by Theorem \ref{notin}, there exist
$P\in  \cP_r$ and $c\in\IC$ with $|c| = 2\|A\|$
such that $\lambda\notin\Lambda_{k-r}(A+cP)$.
\qed
\fi

\iffalse
\it Proof. \rm The inclusion ($\subseteq$) follows from Theorem
\ref{inclusion}.
Suppose $\lambda\notin\Lambda_{k}(A)$. If $\lambda\notin \Omega_k(A)$, then
there exists $t\in \IR$ such that
$\Re\(e^{it}\lambda\)>\lambda_k(\Re\(e^{it}A\))$.
Let $e^{it}A=H+iG$ be with $H = H^*$ and $G=G^*$.
Then $H$ has an operator matrix
$\diag(\lambda_1, \dots, \lambda_m) \oplus H_2$
with $m \le k-1$
such that $\sup \sigma(H_2) < \Re(e^{it}\lambda)$.
Let $F= -2e^{-it}\|A\|(I_m \oplus [I_{k-m} \oplus 0])\in \BH$.
Then  $\Re(e^{it}(A+F))$ has fewer than $k-r$
eigenvalues larger than or equal to $\Re(e^{it}\lambda)$.
Hence,
$\lambda \notin \Omega_{k-r}(A+F)\Ra\lambda \notin  \Lambda_{k-r}(A+F)$.

If $\lambda\in \Omega_k(A)$, then by Theorem \ref{notin}, there exist
$P\in  \cP_r$ and $c\in\IC$ with $|c| = 2\|A\|$
such that $\lambda\notin\Lambda_{k-r}(A+cP)$.
\qed
\fi

\medskip
If $A\in\BH$ is self-adjoint, Theorem \ref{notin} reduces to the following corollary.

\begin{corollary}\label{adjoint}
Suppose $A\in \BH$ is self-adjoint and $1\le r \le k$.
Then $\lambda\in \Lambda_k(A)$ if and only if
$A$ can be decomposed into $\tilde A_1 \oplus \tilde A_2$
such that $\dim \tilde A_1 < k$, $W(\tilde A_1) \subseteq L$
and $W(\tilde A_2) \subseteq \IR \setminus L$,
where $L = [\lambda, \infty)$ or $(-\infty,\lambda]$.
\end{corollary}

Using a similar argument as in the proof of Theorem \ref{normal},
an analogue result can also be obtained for self-adjoint operators.

\begin{theorem}\label{adjoint2}
Suppose $A\in \BH$ is self-adjoint and  $1 \le r < k\le \dim \cH$.
Let $\cS$ be a subset of $\cF_r$ containing the set
\ $\{ \pm 2\|A\|P: P\in \cP_r \}.$
 Then
$$\Lambda_k(A) = \cap\{\Lambda_{k-r}(A+F): F \in \cS\}.$$
\end{theorem}

\medskip
In \cite[Proposition 2.3]{R}, the author showed that
$\Lambda_k(A) \subseteq \cap_{X \in \cV_{k-1}} W(X^*AX)$,
where $\cV_m$ is the set of $X: \cH \rightarrow \cH$
with $X^*X = I_{\cH}$ and $X(\cH)=\cH_1^\perp$ for some subspace $\cH_1$ of $\cH$
satisfying $\dim \cH_1 \le m$. In general, we have the following.

\begin{proposition}
Suppose $A \in \BH$ and $1\le r < k<\infty$. Then
$$\Lambda_k(A) \subseteq \cap\{\Lambda_{k-r} (X^*AX): X\in \cV_r\}.$$
\end{proposition}

\it Proof. \rm Let $\lambda \in \Lambda_k(A)$. Then there exists a rank  $k$ orthogonal projection $P$ such that $PAP=\lambda P$. Suppose  $X\in \cV_r$. Then there exists a subspace $\cH_1$ of $\cH$ with 
 $\dim \cH_1 \le r$ satisfying $X^*X = I_{\cH}$ and $X(\cH)=\cH_1^\perp$. Therefore, dim $\(P(\cH)\cap \cH_1^\perp\)\ge k-r$. Choose a $k-r$ dimensional subspace $\cH_2$ of $\cH$ such that $X\(\cH_2\)\subseteq  P(\cH)\cap \cH_1^\perp$. Let $\{y_i\}_{i=1}^{k-r}$ be an orthogonal basis of $\cH_2$. Then $\{X\(y_i\)\}_{i=1}^{k-r}$ is an orthonormal subset of $P(\cH)$. So,  for $1\le i,\ j\le k-r$, we have
 $$\la X^*AX y_i,y_j\ra =\la A\(X y_i\),\(Xy_j\)\ra =\delta_{i\,j}\lambda.$$
Hence, $\lambda\in \Lambda_{k-r} (X^*AX)$.\qed

Using Theorems \ref{finite} and \ref{normal}, we have 

\begin{corollary}
Suppose $A \in \BH$ and $1\le r < k<\infty$.
If $\dim \cH < \infty$ or $A$ is normal, then
$$\Lambda_k(A) = \cap\{\Lambda_{k-r} (X^*AX): X\in \cV_r\}.$$
\end{corollary}
\it Proof. \rm
For each $F \in \cF_r$, there is $X\in \cV_r$ such that $X^*FX = 0$.
Then
$$\Lambda_k(A) \subseteq \bigcap_{X\in \cV_r} \Lambda_{k-r}(X^*AX)
\subseteq \bigcap_{F\in \cF_r} \Lambda_{k-r}(X^*(A+F)X)
\subseteq \bigcap_{F\in \cF_r} \Lambda_{k-r}(A+F).$$
By Theorems \ref{finite} and \ref{normal}, the inclusions are indeed equalities.
\qed

Similarly, using Theorems \ref{interior} and Corollary \ref{closure}, we have the last corollary in this section.
\begin{corollary}
Suppose $A \in \BH$ and $1\le r < k<\infty$. Then
$$\Omega_k(A) = \cap\{\Omega_{k-r} (X^*AX): X\in \cV_r\}.$$
\end{corollary}

\iffalse
\begin{remark} The results in this section holds for any $\cS$ satisfying
$$\{2e^{i\xi}\|A\|P: P\in \cP_r \mbox{ and }\xi \in [0, 2\pi) \}
\subseteq  \cS\subseteq \cF_r $$
\end{remark}
\fi

\section{Results on $\Lambda_\infty(A)$}
\setcounter{equation}{0}

Suppose $\cH$ is infinite dimensional and $A\in\BH$. It is clear
that $\Lambda_\infty(A)$ can be viewed
as the set of $\lambda\in \IC$ for which there
exists an infinite orthonormal set $\{x_i\in\cH:i\ge 1\}$ such that
$\la Ax_i,x_j\ra =\delta_{i\,j}\lambda$ for all $i,j\ge 1$.
Extend the definition of $\Omega_k(A)$ to
$$\Omega_\infty(A) = \bigcap_{\xi \in [0, 2\pi)}
\left\{\mu \in \IC: \Re(e^{i\xi}\mu) \le
\lambda_k(\Re(e^{i\xi} A)) \mbox{ for all }k\ge 1\right\}.$$

We have the following result.

\begin{theorem}\label{inf}  Suppose dim $\cH$ is infinite and $A\in \BH$.
Let $\cS$ be a set of finite rank operators on $\BH$
containing the set
$$\{2 e^{i\xi}\|A\|P: \xi \in [0, 2\pi), \
P \mbox{ is a finite rank orthogonal projection}\}.$$
Then we have the following equalities.
\begin{enumerate}
\item
$\Omega_{\infty}(A) =\bigcap_{k\ge 1}\Omega_k(A)
= \bigcap \{\cl({W(A+F)}): F \in \cS \}=W_e(A)$. 
\item $\Lambda_{\infty}(A) =\bigcap_{k\ge 1}\Lambda_k(A)
= \bigcap \{W(A+F): F \in \cS\}$.
\end{enumerate}
\end{theorem}

\it Proof. \rm (1) By the definition of $\Omega_\infty(A)$, we have
$\Omega_{\infty}(A) =\bigcap_{k\ge 1}\Omega_k(A).$
By (\ref{ess}) and Corollary \ref{closure}, we have
\begin{eqnarray*}
W_e(A) & = &
\bigcap\{\cl({W(A+F)}): F \in \BH \hbox{ has finite rank }\} \\
& = &
\bigcap_{k\ge 1} \bigcap \{\cl({W(A+F)}): F \in \BH \hbox{ has rank }
k-1\} \\
&=& \bigcap_{k\ge 1} \Omega_k(A)\\
&=& \bigcap_{k\ge 1}
\bigcap \{\cl({W(A+F)}): F \in \cS \hbox{ has rank } k-1\} \\
& = & \bigcap \{\cl({W(A+F)}): F \in \cS \}.
\end{eqnarray*}
So, the second and third equalities in (1) hold.

\iffalse
The first equality follows from definition..

Suppose $\lambda \in  \bigcap_{k\ge 1}\Omega_k(A)$. Let $F\in \BH$
have finite rank $r$. Then for every $\xi\in\IR$,
since $\lambda\in \Omega_{2r+1}(A)$, we have
$$\Re(e^{i\xi}\lambda)\le \lambda_{2r+1}(\Re (e^{i\xi}A))
\le \lambda_{1}( \Re(e^{i\xi}(A+F))).$$
Therefore, $\lambda\in \cl({W(A+F)})$.
This proves that
$$\bigcap_{k\ge 1}\Omega_k(A)\subseteq \bigcap
\{\cl({W(A+F)}): F \in \BH \hbox{ has finite rank} \}.$$

Now, suppose $\lambda\notin  \bigcap_{k\ge 1}\Omega_k(A)$.
Then $\lambda\notin \Omega_k(A)$ for some $k\ge 1$.
So there exists $\xi\in\IR$ such that
$\lambda_k(\Re(e^{i\xi}A))<\Re(e^{i\xi}\lambda)$.
Let $B= \Re(e^{i\xi}A)$ and  denote $\lambda_j(B)$ by
$\lambda_j$, $1\le j\le k$.. There exists $r\le k-1$ such that
$\lambda_1\ge \cdots \ge\lambda_r$ are eigenvalues of $B$ with
$\lambda_r\ge \lambda>\lambda_{r+1}$. Choose a set of orthonormal vectors
$\{u_j:1\le j\le r\}$ such that $Bu_j=\lambda_j u_j$ for  $1\le j\le r$.
Define $F\in \BH$ by
$Fx=-e^{-i\xi} 2\|A\| \(\sum_{j=1}^r \la x,u_j\ra u_j\)$.
Then $F \in \cS$ and
$$\lambda_1\(\Re(e^{i\xi}(A+F))\)
=\lambda_{r+1}(\Re(e^{i\xi}A))<\lambda.$$
Therefore, $\lambda\not \in \cl({W(A+F)})$.
This proves
$$\bigcap \{\cl({W(A+F)}): F \in \BH \hbox{ has finite rank}\}
\subseteq W_e(A).$$
\fi

(2)  By Theorem 4 in \cite{AS},we have
\begin{equation}\label{4.5}
\Lambda_{\infty}(A) =
\bigcap\{W(A+F): F \in \BH \hbox{ has finite rank } \}\,.
\end{equation}
Clearly, we have the inclusion
$$\Lambda_{\infty}(A) \subseteq \bigcap_{k\ge 1}\Lambda_k(A).$$
To prove the reverse inclusion,
suppose $\lambda\in\bigcap_{k\ge 1}\Lambda_k(A)$.
Let  $F\in \BH$  of rank $m$. Choose $k\ge m+ 1$. Then $\lambda\in
\Lambda_k(A)$. By Theorem \ref{inclusion}, we have
$$\lambda\in \Lambda_k(A)\subseteq
\Lambda_1(A+F) = W(A+F).$$
Hence,
$$\bigcap_{k\ge 1}\Lambda_k(A)\subseteq
\bigcap\{W(A+F): F \in \BH \mbox{ is  of finite rank} \}
= \Lambda_\infty(A).$$
Thus, we get the first equality in (2).

Next, we show that one only needs to use $F \in \cS$ for the intersection
on the right side of (\ref{4.5}). To this end, note that
$$\bigcap\{W(A+F): F \in \BH \mbox{ is  of finite rank} \}
\subseteq \bigcap\{W(A+F): F \in \cS\}.$$
To prove the reverse inclusion, assume that
$$\lambda \notin
\bigcap\{W(A+F): F \in \BH \mbox{ is  of finite rank} \}.$$
If $\lambda \notin \Omega_\infty(A)$, then there is a finite rank
$F \in \BH$ such that $\lambda \notin \cl(W(A+F))$
and hence  $\lambda \notin W(A+F)$.
So, assume that
\begin{equation}\label{4.2}
\lambda \in W_e(A) \subseteq W(A) \hbox{ and thus }
|\lambda| \le \sup\{ |\mu|: \mu \in W(A)\} \le \|A\|.
\end{equation}
Then there is $\xi \in [0, 2\pi)$ and a finite rank operator
$F \in \BH$ such that
\begin{equation}\label{cknew}
e^{i\xi}W(A+F-\lambda I) \subseteq \{\mu\in \IC: \im(\mu) < 0\} \cup \cL
\end{equation}
with $\cL = (0, \infty)$ or $\cL = (-\infty,0)$.
We may replace $A$ by $e^{i\xi}A$ and assume that $\xi = 0$.
Without loss of generality, assume that $\cL = (-\infty,0)$.

Let $\lambda = a+ib$ with $a, b\in \IR$
and $A = H+iG$ with $H = H^*$ and $G = G^*$.
Since (\ref{cknew}) holds
with $\xi = 0$ for a finite rank operator $F \in \BH$,
there is $r$ not larger than the rank of $\im F$ such that
$G$ has an operator matrix of the form
\begin{equation} \label{Gform}
\diag(g_1, \dots, g_r) \oplus bI_s \oplus G_2
\end{equation}
with
$g_1 \ge \cdots \ge g_r > b$, $W(G_2) \subseteq (-\infty,b)$
and $0 \le s \le \infty$.
By (\ref{4.2}), we have
$$g_1 -b \le |g_1|+|b| \le 2\|G\| \le 2\|A\|.$$
We consider two cases.

\bf Case 1. \rm
Suppose $g_1 - b = 2\|A\|$. Then $g_1 = \|A\| = -b$.
Since $\lambda\in W_e(A)$ and
$$\|A\| = |b| \le |a+ib| = |\lambda| \le \|A\|,$$
it follows that
$$a = 0 \quad  \hbox{ and } \quad \lambda = ib = -i\|A\|$$
is the only element in $\cl(W(A)) \cap \{\mu \in \IC: \im(\mu) \le -\|A\|\}$.
Thus, $G_2$ in (\ref{Gform}) is vacuous, i.e., $G$
has operator matrix $\diag(g_1, \dots, g_r) \oplus bI_s$.
Using the same basis, we let $H$ have the operator matrix
$$\begin{bmatrix}H_{11}& H_{12} \cr H_{12}^* & H_{22}\end{bmatrix}.$$
Since $\|H_{22}+ibI\| \le \|A\| = |b|$, we see that
$H_{22} = 0$. By the fact that
$$|b|^2 = \|A^*A\| = \|(H+iG)^*(H+iG)\|,$$
we see that $H_{12}$ is zero as well.
Thus, $A$ has operator matrix
$$A_1 \oplus ibI_s \qquad \hbox{ with } \ A_1 \in M_m.$$
Since (\ref{cknew}) holds for a finite operator $F$
with $\xi = 0$ and $\cL = (-\infty,0)$,
we see that $s \ne \infty$.
But then $\dim \cH$ is finite, which is a contradiction.

\bf Case 2. \rm  Suppose $g_1 - b < 2 \|A\|$.
If $s$ is finite in (\ref{Gform}), then
$$\tilde F = i2\|A\|(I_{r+s} \oplus 0) \in \cS
\quad \hbox{ and } \quad
W(A- \tilde F) \subseteq \{\mu \in \IC: \im(\mu) < b\}.$$
Thus,  $\lambda = a+ib \notin W(A-\tilde F)$.

Next, assume that $s = \infty$.
Suppose the compression of $H$ on the null space of
$G-bI$ equals $H_0$. Then there is a positive integer $m$ such that
$H_0$ has operator matrix
$\diag(h_1, \dots, h_m) \oplus H_1	$
such that $h_1 \ge \cdots \ge h_m \ge a$ and
$W(H_1	) \subseteq (-\infty,a)$.
Otherwise, (\ref{cknew}) cannot hold for a finite operator $F$
with $\xi = 0$ and $\cL = (-\infty,0)$.
Let $\tilde F = i2\|A\|(I_{r+m} \oplus 0) \in \cS$,
and let $\hat A = A - \tilde F - \lambda I$.
Then $\im(\hat A) = (\hat A - \hat A^*)/2i$
has an operator matrix $\hat G_1 \oplus 0_{s-m} \oplus \hat G_2$
with $W(\hat G_1 \oplus \hat G_2) \subseteq (-\infty,0)$.
Moreover, the compression of $\re(\hat A) = (\hat A + \hat A^*)/2$ on the null space of
$\im(\hat A)$ equal $H_1 - aI$.
As a result,
if  $\mu = \langle \hat A x,x\rangle \in W(\hat A)$
has imaginary part 0, then $x$ must lie in the null space of
$\im(\hat A)$, and hence the real part of $\mu$ lies
in $W(H_1 - aI) \subseteq (-\infty,0)$.
Thus, $0 \notin W(\hat A)$, equivalently,
$\lambda \notin W(A - \tilde F)$.
Consequently,
$$\bigcap\{W(A+F): F \in \cS\} \subseteq
\bigcap\{W(A+F): F \in \BH \mbox{ is  of finite rank} \}.$$
\vskip -.3in\qed

\medskip

In \cite{R}, Martinex-Avendano asked whether
$\Lambda_\infty(A) = \bigcap_{k\ge 1} \Lambda_k(A)$.
Assertion (2) answers the question affirmatively.

\begin{theorem} \label{last} Suppose $A \in \BH$, where $\cH$ is infinite
dimensional.
Then
$$\int(\Omega_\infty(A)) \subseteq
\Lambda_\infty(A) \subseteq \Omega_\infty(A).$$
Moreover, $\cl(\Lambda_\infty(A)) = \Omega_\infty(A)$ if and only if
$\Lambda_\infty(A) \ne \emptyset$.
\end{theorem}

\it Proof. \rm
By the Corollary after Theorem 4 in \cite{AS},
we see that   $\int(\Omega_\infty(A)) \subseteq \Lambda_\infty(A)$.
The inclusion $\Lambda_\infty(A) \subseteq \Omega_\infty(A)$ is clear.

Note that $\Omega_\infty(A)$ is always a non-empty compact convex set.
If $\Lambda_\infty(A) = \emptyset$ then
$\cl(\Lambda_\infty(A)) \ne \Omega_\infty(A)$.
Conversely, suppose $\Lambda_\infty(A) \ne \emptyset$.
If $\int(\Lambda_\infty(A))=\int(\Omega_\infty(A))$ is non-empty,
then $\cl(\Lambda_\infty(A)) = \Omega_\infty(A)$.
If $\int(\Omega_\infty(A))$ is empty, then
$\Omega_\infty(A) = \{\mu \}$ is a singleton
and so is the non-empty set $\Lambda_\infty(A)$. Hence
$\cl(\Lambda_\infty(A))=\Lambda_\infty(A) = \{\mu\}$.
\qed

The next example show that $\Lambda_\infty(A)$ may indeed be empty
so that $\cl(\Lambda_\infty(A)) \ne \Omega_\infty(A)$.

\begin{example}\label{compact}
Let $A=\bigoplus_{n\ge 2}\diag\({e^{i\pi/n}}/{n},-{1}/{n}\) \in \BH$.
Then $\Omega_\infty(A) = \{0\}$ but
$0\notin \Lambda_1(A)$ so that $\Lambda_\infty (A)= \cap
\{\Lambda_k(A): k = 1, 2, \dots \} = \emptyset$.
On the other hand, if $B = A \oplus 0_\cH$, then $\Lambda_\infty(B) = \{0\}$.
\end{example}

From the proof of Theorem \ref{last}, we see that
if $\Lambda_\infty(A)$ is a singleton, then $\Omega_\infty(A)$
is also a singleton, which can happen if and only
if $A-\mu I$ is a compact operator for some $\mu \in \IC$
by the corollary after Lemma 3 in \cite{AS}.
In connection to this comment and Example \ref{compact},
we have the following.

\begin{proposition}
Let $A \in \BH$ and $\mu \in \IC$ be such that
$A-\mu I$ is  compact. Then the following are equivalent.
\begin{itemize}
\item[(a)] $\Lambda_\infty(A)$ is non-empty.
\item[(b)] $\Lambda_\infty(A) = \{\mu\}$.
\item[(c)] $\mu \in \Lambda_k(A)$ for each $k = 1,2,\dots$.
\end{itemize}
\end{proposition}

\it Proof. \rm
The implications ``(a) $\iff$ (b)'' is clear.
We have ``(c) $\iff$ (b)''
because $\Lambda_\infty(A) = \bigcap_k \Lambda_k(A)$ by Theorem \ref{inf}.
\qed

\medskip\noindent
{\bf Acknowledgment}

Li is an honorary professor of the University of Hong
Kong. His research was partially supported by an USA NSF grant
and a HK RGC grant. This research was done while the three authors
were visiting the University of Hong Kong in the summer of
2007 supported by a HK RGC grant.
They would like to thank the staff of the Mathematics Department
for their hospitality.

\bibliographystyle{amsplain}


\end{document}